\documentclass[12pt,leqno]{article}
\usepackage{graphicx}
\usepackage{hyperref}
\usepackage{amssymb,amsmath,amsthm}
\usepackage{epstopdf}
\usepackage{mathrsfs}
\usepackage{xypic}
\usepackage{comment}
\usepackage{mathabx,epsfig}
\usepackage{comment}
\usepackage{parskip}
\usepackage{cancel}

\def\RR{\mathbb{R}}

\def\p{\mathbf{p}}

\def\r{\mathbf{r}}
\def\s{\mathbf{s}}
\def\t{\mathbf{t}}
\def\u{\mathbf{u}}
\def\v{\mathbf{v}}
\def\x{\mathbf{x}}

\def\0{\mathbf{0}}

\allowdisplaybreaks

\makeatletter
\renewcommand\@biblabel[1]{}
\makeatother

\title{Where is the cone?}
\author{Drew Armstrong \\ \url{www.math.miami.edu/~armstrong}}
\date{August 23, 2017}

\begin{document}

\maketitle

Real quadric curves are often referred to as ``conic sections," implying that they can be realized as plane sections of circular cones. However, it seems that the details of this equivalence have been partially forgotten by the mathematical community. The definitive analytic treatment was given by Otto Staude in the 1880s and a non-technical description was given in the first chapter of Hilbert and Cohn-Vossen's {\em Geometry and the Imagination} (1932). In this note we will prove an elegant theorem that completely answers the question for nondegenerate quadric curves (ellipses and hyperbolas). The theorem is easy to state but it was very difficult to track down, and we still have not seen a proof in the literature. Our hope is to revive the lost knowledge of ``conic sections" by providing the slickest possible modern treatment, by using standard linear algebra that was not standard in 1932.

\tableofcontents

\section{Introduction}
\label{sec:introduction}

A {\em real quadric curve} is defined by an equation of the form
\begin{equation*}
ax^2+bxy+cy^2+dx+ey+f=0
\end{equation*}
for some real numbers $a,b,c,d,e,f\in\RR$. One often sees real quadric curves described as ``conic sections." This terminology suggests that we can view this curve as the intersection of the $x,y$-plane with a circular cone living in $x,y,z$-space. However, it is very rare to see the details of this spelled out. The purpose of the present note is to answer the question: {\bf where is the cone?} More specifically, we will answer the following three questions:
\begin{itemize}
\item From which points in space does a given quadric curve look like a circle?\\ (This is the apex of the cone.)
\item In which direction should we look to see this circle?\\(This is the axis of symmetry of the cone.)
\item How big is the circle?\\(This is the angle of aperture of the cone.)
\end{itemize}

The answers to these questions were well known to geometers in the late nineteenth and early twentieth centuries. However, based on my several years of internet searching, it seems that the answers are not well known today. Eventually I found a clue in footnote 4 on page 24 of Hilbert and Cohn-Vossen's {\em Geometry and the Imagination} (1932).\footnote{For the precise statement see Section \ref{sec:symmetries} below.} This led me to a genre of textbooks on the ``analytic geometry of three dimensions" that were written around the same time. I found the text of D.M.Y.~Sommerville (1934) particularly helpful.

Without further ado, here is the main theorem describing cones over quadric curves.

\bigskip

{\bf Main Theorem.}\footnote{As stated this result applies only to {\em non-degenerate} and {\em central} quadric curves (ellipses and hyperbolas). Analogous results for non-central quadric curves (parabolas) and degenerate curves (line pairs) can be obtained from limiting arguments that are not so interesting, so I will omit them.} Let $a,b,c\in\RR$ be distinct real numbers and consider the following space curves in the three principal coordinate planes:
\begin{align*}
x^2/(a-c)+y^2/(b-c)=1 \quad\text{and}\quad z=0,\\
x^2/(a-b)+z^2/(c-b)=1 \quad\text{and}\quad y=0,\\
y^2/(b-a)+z^2/(c-a)=1 \quad\text{and}\quad x=0.
\end{align*}
In general, two of these curves are real and the third is imaginary. If $\u$ is any point on one of the real curves then the cone from $\u$ to the other real curve is circular. Furthermore, the axis of symmetry of this cone is the tangent line to the first curve at $\u$.\footnote{and from this one can easily compute the angle of aperture} The cone from any {\bf other} point in space to any one of the real curves is {\bf not} circular. \hfill ///

\bigskip

\begin{figure}
\begin{center}
\includegraphics[scale=0.5]{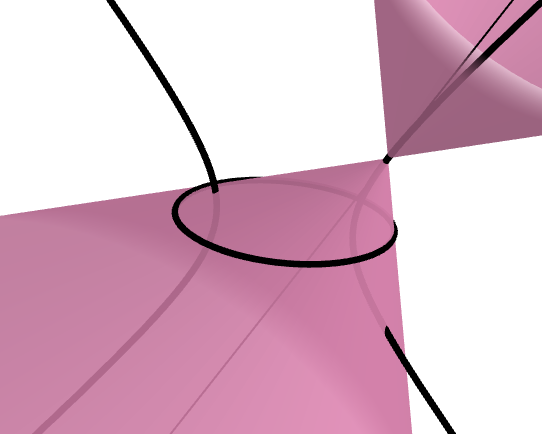}
\end{center}
\caption{Realizing an ellipse as a conic section.}
\label{fig:ellipse}
\end{figure}

\begin{figure}
\begin{center}
\includegraphics[scale=0.5]{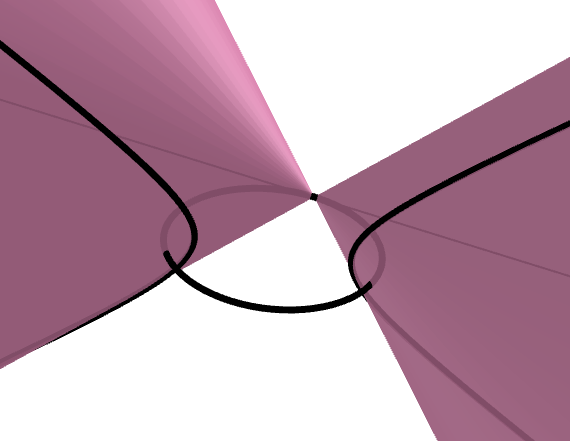}
\end{center}
\caption{Realizing a hyperbola as a conic section.}
\label{fig:hyperbola}
\end{figure}

For example, if $c<b<a$ then we have an ellipse in the $x,y$-plane and a hyperbola in the $x,z$-plane. Figure \ref{fig:ellipse} shows a typical circular cone over the ellipse and Figure \ref{fig:hyperbola} shows a typical circular cone over the hyperbola.\footnote{These pictures were made with GeoGebra.}

This theorem provides a completely satisfying answer to the question ``where is the cone?" and I am surprised by how difficult it was to track down. After reading the result in Hilbert and Cohn-Vossen I found it stated without proof in several textbooks of the period.\footnote{It even appears as Exercise 6 on page 100 of Barry Spain's {\em Analytic Quadrics} (1960).} More recently, the result appears as Theorem 4.2.1 in {\em The Universe of Conics} (2016), where the authors provide a case-by-case synthetic proof. However, I still have not found an analytic proof written down anywhere.

My goal in this note is to provide the slickest possible analytic treatment for the statement in Hilbert and Cohn-Vossen's footnote by employing standard linear algebra that was not standard in 1932. I don't claim to have found the ``book proof" but hopefully I have done something useful to fill this surprising gap in the literature.

\section{The Equation of a Circular Cone}
\label{sec:circcone}

To give an analytic treatment of conic sections we must first give an analytic treatment of cones. A {\em circular cone} in $\RR^3$ is specified by the following data:
\begin{itemize}
\item A point in space $\u\in\RR^3$ defining the apex,
\item A unit vector $\r\in\RR^3$ defining the axis of symmetry,
\item An angle $\theta\in[0,\pi/2]$ defining the aperture of the cone.\footnote{The extreme values $\theta=0$ and $\theta=\pi/2$ correspond to a line and a plane, respectively.}
\end{itemize}
Geometrically, the cone consists of all points $\x\in\RR^3$ such that the line connecting $\x$ and $\u$ makes an angle of $\theta$ (or $\pi-\theta$) with the axis of symmetry. Algebraically, we can express this situation with the dot product. Since $\r$ is a unit vector we have
\begin{equation*}
(\x-\u)^T\r=\|\x-\u\|\cdot\|\r\|\cdot\cos\theta=\|\x-\u\|\cdot\cos\theta. \tag{$\ast$}
\end{equation*}
Furthermore, since every $1\times 1$ matrix is symmetric we have
\begin{equation*}
(\x-\u)^T\r=\left[(\x-\u)^T\r\right]^T=\r^T(\x-\u).
\end{equation*}
Then by squaring both sides of ($\ast$) we obtain the following equation for the circular cone:
\begin{align*}
[(\x-\u)^T\r]^2&=\|\x-\u\|^2\cdot\cos^2\theta \\
[(\x-\u)^T\r]\cdot [\r^T(\x-\u)] &= (\x-\u)^T(\x-\u)\cdot\cos^2\theta \\
(\x-\u)^T\left(\r\r^T\right)(\x-\u)&= (\x-\u)^T\left((\cos^2\theta)I\right)(\x-\u)\\
(\x-\u)^T\left(\r\r^T-(\cos^2\theta)I\right)(\x-\u) &=0.
\end{align*}
In summary, for any unit vector $\r=(r,s,t)$ and angle $\theta$ we define the matrix
\begin{equation*}
C_{\r,\theta}:=\r\r^T-(\cos^2\theta)I=\begin{pmatrix} r^2-\cos^2\theta&rs&rt\\ rs&s^2-\cos^2\theta&st \\ rt& st& t^2-\cos^2\theta\end{pmatrix}.
\end{equation*}
Then the circular cone with apex $\u$, axis $\r$ and aperture $\theta$ is defined by the equation
\begin{equation*}
\boxed{(\x-\u)^T C_{\r,\theta}(\x-\u)=0.}
\end{equation*}

Let us make a few observations about the {\em cone matrix} $C_{\r,\theta}$. Since $\r^T\r=\|\r\|^2=1$ we observe that $\r$ is an eigenvector of $C_{\r,\theta}$ with eigenvalue $1-\cos^2\theta$:
\begin{align*}
C_{\r,\theta}\cdot \r=\left(\r\r^T-(\cos^2\theta)I\right)\r=\r(\r^T\r)-(\cos^2\theta)\r=(1-\cos^2\theta)\r.
\end{align*}
And if $\r'$ is any vector perpendicular to $\r$ (i.e., if $\r^T\r'=0$) then we observe that $\r'$ is an eigenvector of $C_{\r,\theta}$ with eigenvalue $-\cos^2\theta$:
\begin{align*}
C_{\r,\theta}\cdot \r'=\left(\r\r^T-(\cos^2\theta)I\right)\r'=\r(\r^T\r')-(\cos^2\theta)\r'=-(\cos^2\theta)\r'.
\end{align*}
We conclude that the cone matrix $C_{\r,\theta}$ has eigenvalues
\begin{equation*}
1-\cos^2\theta,\quad -\cos^2\theta,\quad -\cos^2\theta.
\end{equation*}
Conversely, let $C^T=C$ be {\bf any} real symmetric matrix with the same eigenvalues. In this case I claim that $C=C_{\r,\theta}$ for some unit vector $\r\in\RR^3$. Indeed, since any real symmetric matrix is orthogonally diagonalizable (by the Principal Axis Theorem), we know that there exists a real orthogonal matrix $Q^T=Q^{-1}$ such that
\begin{equation*}
Q^T C Q = \begin{pmatrix} 1-\cos^2\theta &0&0 \\ 0&-\cos^2\theta & 0 \\ 0&0&-\cos^2\theta\end{pmatrix}=\begin{pmatrix} 1&0&0\\0&0&0\\0&0&0\end{pmatrix}-(\cos^2\theta)I.
\end{equation*}
Now let $\r$ be the first column of the matrix $Q$, i.e., the eigenvector of $C$ corresponding to eigenvalue $1-\cos^2\theta$. Since $QQ^T=I$ we see that $\r$ is a unit vector and we observe that
\begin{equation*}
C=Q\left[\begin{pmatrix} 1&0&0\\0&0&0\\0&0&0\end{pmatrix}-(\cos^2\theta)I\right]Q^T=Q\begin{pmatrix} 1&0&0\\0&0&0\\0&0&0\end{pmatrix}Q^T-(\cos^2\theta)QQ^T=\r\r^T-(\cos^2\theta)I
\end{equation*}
as desired.

\section{Quadric Cones in General}
\label{sec:quadcone}

More generally, for any real symmetric matrix $C^T=C$ we consider the quadratic equation
\begin{equation*}
(\x-\u)^TC(\x-\u)=0.
\end{equation*}
Clearly the point $\x=\u$ is a solution. If the (necessarily real) eigenvalues of $C$ are either all negative or all positive then $\x=\u$ is the {\bf only (real) solution}. We are interested in the case when the eigenvalues of $C$ are nonzero and not all of the same sign. In this case we say that $C$ is {\em indefinite} and we say that the equation $(\x-\u)^TC(\x-\u)=0$ defines a {\em quadric cone}.

The following theorem says that the reflection symmetries of a quadric cone are the same as the eigenvectors of its matrix. In fact it is true that {\bf every} symmetry of the cone is a product of {\bf reflection} symmetries,\footnote{This is the famous Cartan-Dieudonn\'e Theorem.} but we won't need that result.

\bigskip

{\bf Theorem.} Let $C^T=C$ be an indefinite real symmetric matrix (i.e., whose eigenvalues are nonzero and not all of the same sign). For any vector $\p\in\RR^3$ the following are equivalent:
\begin{itemize}
\item We have $C\p=c\p$ for some scalar $c\in\RR$.
\item The quadric cone $(\x-\u)^TC(\x-\u)=0$ is symmetric about the plane $\u+\p^\perp$ which is perpendicular to $\p$ and passes through the apex $\u$.
\end{itemize}
\hfill ///

\bigskip

To prove this we require a general lemma on quadratic forms, which is an easy consequence of Hilbert's Nullstellensatz. To keep the treatment self-contained I will present an elementary proof that I learned from John H.~Elton (2009).

\bigskip

{\bf Lemma.} Let $A^T=A$ and $B^T=B$ be real symmetric matrices of the same size, with $A$ indefinite. If the cone of $A$ is contained in the cone of $B$, i.e., if for all vectors $\x$ we have
\begin{equation*}
\x^TA\x=0 \quad\Longrightarrow\quad \x^TB\x=0,
\end{equation*}
then it follows that $B=\lambda A$ for some real scalar $\lambda\neq 0$.\hfill ///

\bigskip

{\bf Proof of the Lemma.} We will prove the result for $3\times 3$ matrices. (See Elton for the general case.) By the Principal Axis Theorem there exists an orthogonal matrix $Q^T=Q^{-1}$ such that $Q^TAQ$ is diagonal:
\begin{equation*}
Q^TAQ=\begin{pmatrix} \lambda_1 & & \\ & \lambda_2 & \\ &&\lambda_3\end{pmatrix}.
\end{equation*}
Since $Q$ is invertible, it is enough to show that $Q^TBQ$ is a scalar multiple of $Q^TAQ$. Furthermore, since $A$ is indefinite we can assume without loss of generality (replacing $A$ by $-A$ if necessary) that the eigenvalues of $A$ satisfy $\lambda_1 > 0 > \lambda_2,\lambda_3$. Let us define the positive real numbers
\begin{equation*}
\ell_1:=1/\sqrt{\lambda_1}, \quad \ell_2:=1/\sqrt{-\lambda_2}\quad\text{and}\quad \ell_3:=1/\sqrt{-\lambda_3}.
\end{equation*}
We observe that the following two equations hold:
\begin{equation*}
0=\begin{pmatrix} \ell_1 & \pm \ell_2 & 0\end{pmatrix} \begin{pmatrix} \lambda_1 & & \\ & \lambda_2 & \\ &&\lambda_3\end{pmatrix} \begin{pmatrix} \ell_1\\ \pm \ell_2 \\ 0 \end{pmatrix} = \begin{pmatrix} \ell_1 & \pm \ell_2 & 0\end{pmatrix}Q^TAQ \begin{pmatrix} \ell_1\\ \pm \ell_2 \\ 0 \end{pmatrix}.
\end{equation*}
Then the implication $(\x^TA\x=0)\Rightarrow(\x^TB\x=0)$ with $\x^T=\begin{pmatrix} \ell_1 &\pm\ell_2&0\end{pmatrix} Q^T$ tells us that
\begin{equation*}
0=\begin{pmatrix} \ell_1 & \pm \ell_2 & 0\end{pmatrix}Q^TBQ \begin{pmatrix} \ell_1\\ \pm \ell_2 \\ 0 \end{pmatrix}.
\end{equation*}
By writing $b_{ij}$ for the $(i,j)$-the entry of $Q^TBQ$ these two equations become\begin{equation*}
\left\{
\begin{array}{ccccccc}
b_{11}/\lambda_1 & - & b_{22}/\lambda_2 & + & 2b_{12}\ell_1\ell_2 & = & 0, \\
b_{11}/\lambda_1 & - & b_{22}/\lambda_2 & - & 2b_{12}\ell_1\ell_2 & = & 0.
\end{array}
\right.
\end{equation*}
By adding these equations we find that $b_{22}=\lambda_2(b_{11}/\lambda_1)$ and hence $b_{12}=0$. Then a similar argument shows that
\begin{equation*}
0=\begin{pmatrix} \ell_1 & 0 & \pm\ell_3\end{pmatrix}Q^TBQ \begin{pmatrix} \ell_1\\  0 \\ \pm \ell_3 \end{pmatrix},
\end{equation*}
which implies that $b_{33}=\lambda_3(b_{11}/\lambda_1)$ and $b_{13}=0$. At this point we know that
\begin{equation*}
Q^TBQ-\left(\frac{b_{11}}{\lambda_1}\right)Q^TAQ = \begin{pmatrix} 0&0&0 \\ 0&0&b_{23} \\ 0&b_{23}&0 \end{pmatrix}
\end{equation*}
and it remains only to show that $b_{23}=0$. To do this we fix any real numbers $(\alpha,\beta,\gamma)$ with the properties $\alpha^2=\beta^2+\gamma^2$ and $\beta\gamma\neq 0$.\footnote{Elton chooses $(\alpha,\beta,\gamma)=(5,4,3)$.} Then we obtain
\begin{equation*}
\begin{pmatrix} \alpha\ell_1&\beta\ell_2&\gamma\ell_3\end{pmatrix} Q^TAQ \begin{pmatrix} \alpha\ell_1\\ \beta\ell_2\\ \gamma\ell_3\end{pmatrix}=\alpha^2-\beta^2-\gamma^2=0 \quad\Rightarrow\quad \begin{pmatrix} \alpha\ell_1&\beta\ell_2&\gamma\ell_3\end{pmatrix} Q^TBQ \begin{pmatrix} \alpha\ell_1\\ \beta\ell_2\\ \gamma\ell_3\end{pmatrix}=0,
\end{equation*}
which implies that
\begin{equation*}
2\beta\gamma\ell_2\ell_3 b_{23}=\begin{pmatrix} \alpha\ell_1&\beta\ell_2&\gamma\ell_3\end{pmatrix} \begin{pmatrix} 0&0&0 \\ 0&0&b_{23} \\ 0&b_{23}&0 \end{pmatrix} \begin{pmatrix} \alpha\ell_1\\ \beta\ell_2\\ \gamma\ell_3\end{pmatrix}=0.
\end{equation*}
Finally, since $\beta\gamma\neq 0$ we conclude that $b_{23}=0$. This completes the proof.\hfill $\qed$

\bigskip

{\bf Proof of the Theorem.} Recall that $C^T=C$ is an indefinite real symmetric matrix. It is sufficient to prove the theorem in the case when $\p$ is a unit vector and the apex is at the origin: $\u=\0$. Since $\p$ is a unit vector (i.e., $\p^T\p=1$) we observe that $P=\p\p^T$ is the matrix that projects orthogonally onto the line $\RR\p$ and that
\begin{equation*}
R=I-2P=I-2\p\p^T
\end{equation*}
is the matrix that reflects orthogonally across the plane $\p^\perp$.

First let us suppose that $\p$ is an eigenvector of $C$, say $C\p=c\p$. To prove that the reflection $R$ leaves the cone invariant we must show for any vector $\x$ that
\begin{equation*}
\x^TC\x=0 \quad\Rightarrow\quad (R\x)^T C(R\x)=\x^T RCR \x =0.
\end{equation*}
By observing that
\begin{align*}
RCR &= (I-2\p\p^T)C(I-2\p\p^T) \\
&= C-2\p(\p^TC) - 2(C\p)\p^T + 4\p\p^T(C\p)\p^T \\
&= C -2\p(c\p^T)-2(c\p)\p^T+4\p\p^T(c\p)\p^T \\
&= C -2c\p\p^T-2c\p\p^T+4c\p(\p^T\p)\p^T \\
&= C -2c\p\p^T-2c\p\p^T+4c\p\p^T \\
&= C,
\end{align*}
we see that the statement is vacuously true. Conversely, let us suppose that
\begin{equation*}
\x^TC\x=0 \quad\Rightarrow\quad \x^T RCR\x =0
\end{equation*}
for all vectors $\x$. Since $C$ is indefinite, the Lemma tells us that $RCR=\lambda C$ for some scalar $\lambda\in\RR$. Since $\det(C)\neq 0$ and $\det(R)=-1$ this implies that
\begin{align*}
\det(\lambda C) &= \det(RCR) \\
\lambda^3 \det(C) &= \det(R)^2 \det(C) \\
\lambda^3 \det(C) &= \det(C) \\
\lambda^3 &=1,
\end{align*}
which implies that $\lambda=1$ since $\lambda$ is real. Thus we conclude that $RCR=C$.  Finally, since $R^2=I$ (indeed, $R$ is a reflection) we observe that
\begin{align*}
RC &= CR \\
(I-2\p\p^T)C &= C(I-2\p\p^T) \\
C-2\p(C\p)^T &= C-2(C\p)\p^T \\
\p(C\p)^T &= (C\p)\p^T.
\end{align*}
Let $\p_i$ and $(C\p)_i$ denote the $i$-th entries of the vectors $\p$ and $C\p$, respectively. By comparing the $(i,j)$-th entry on both sides of the previous equation we find that $\p_i(C\p)_j=(C\p)_i\p_j$ for all $i$ and $j$. Since $\p\neq\0$ there exists some index $i$ such that $\p_i\neq 0$; let us define $c:=(C\p)_i/\p_i$. Then for all indices $j$ we have $(C\p)_j=c\p_j$ and it follows that $C\p=c\p$ as desired. \hfill $\qed$

\bigskip

The following theorem summarizes the results of Sections \ref{sec:circcone} and \ref{sec:quadcone}.

\bigskip

{\bf Theorem.} Let $C^T=C$ be an indefinite real symmetric matrix with (real, nonzero) eigenvalues $\lambda,\mu,\nu$, not all of the same sign. By the Principal Axis Theorem there exists an orthogonal basis of eigenvectors $\r,\s,\t$ corresponding to eigenvalues $\lambda,\mu,\nu$, respectively.
\begin{itemize}
\item If $\lambda,\mu,\nu$ are distinct then the quadric cone $(\x-\u)^TC(\x-\u)=0$ is symmetric with respect to the three mutually perpendicular planes $\u+\r^\perp, \u+\s^\perp, \u+\t^\perp$, called the {\em principal planes} of the cone. The cone is not symmetric with respect to any other plane.
\item If two eigenvalues collide, say $\mu=\nu$, then the cone $(\x-\u)^TC(\x-\u)=0$ is circular with axis of symmetry $\u+\RR\r$ and angle of aperture $\theta$ satisfying
\begin{equation*}
\cos^2\theta=\frac{\mu}{\mu-\lambda}.
\end{equation*}
\end{itemize}
\hfill ///

\bigskip

{\bf Proof.} It remains only to prove the second statement. We assume that the eigenvalues of $C$ are $\lambda,\mu,\mu$ with corresponding orthogonal eigenbasis $\r,\s,\t$. Since the eigenvalues of $\s$ and $\t$ are equal we see that the cone is symmetric with respect to any plane of the form $\u+(a\s+b\t)^\perp$. In other words, the cone has rotational symmetry around the perpendicular axis $\u+\RR\r$. Now let us consider the matrix $C':=C/(\lambda-\mu)$ with eigenvalues
\begin{equation*}
\frac{\lambda}{\lambda-\mu},\quad\frac{\mu}{\lambda-\mu},\quad\frac{\mu}{\lambda-\mu}.
\end{equation*}
Observe that the equations $(\x-\u)^TC(\x-\u)=0$ and $(\x-\u)^T C'(\x-\u)=0$ define the same cone. Since $\lambda$ and $\mu$ have opposite signs we observe that $0<\mu/(\mu-\lambda),\lambda/(\lambda-\mu)<1$ and $\mu/(\mu-\lambda)+\lambda/(\lambda-\mu)=1$, hence there exists a unique angle $\theta\in[0,\pi/2]$ satisfying
\begin{equation*}
\frac{\lambda}{\lambda-\mu}=1-\cos^2\theta \qquad\text{and}\qquad \frac{\mu}{\mu-\lambda}=\cos^2\theta.
\end{equation*}
It follows from the remarks of Section \ref{sec:circcone} that $C'=C_{\r,\theta}$ as desired. \hfill $\qed$

\section{Tangent Cones Over Quadric Surfaces}
\label{sec:tangentcone}

Our goal is to study the cone from a point $\u$ in $x,y,z$-space to a quadric curve lying in the $x,y$-coordinate plane. Surprisingly, it turns out that the best way to do this is to first consider the {\bf tangent cones} from a given point $\u$ to a certain family of {\bf quadric surfaces} in $x,y,z$-space. Then we will allow these quadric surfaces to degenerate to the desired quadric curve.

The quadric surfaces we will consider have the form
\begin{equation*}
\x^T A\x =1,
\end{equation*}
where $A^T=A$ is a real symmetric matrix. If $\u$ is any point, then the {\em tangent cone from $\u$ to the surface} consists of all points $\x$ such that the line $t\x+(1-t)\u$ has double contact with the surface. Observe that the line and surface intersect when
\begin{align*}
(tx+(1-t)\u)^T A (t\x+(1-t)^2\u) &= 1 \\
t^2 \x^TA\x+2t(1-t)\x^tA\u+(1-t)^2\u^TA\u &= 1\\
t^2(\x^TA\x+\u^TA\u-2\x^TA\u)+t(2\x^TA\u-2\u^TA\u)+(\u^TA\u-1) &= 0.
\end{align*}
This quadratic equation in $t$ usually has two distinct roots corresponding to two distinct points of contact with the surface. Double contact occurs when the discriminant vanishes, i.e., when
\begin{align*}
(2\x^TA\u-2\u^TA\u)^2-4(\x^TA\x+\u^TA\u-2\x^TA\u)(\u^TA\u-1) &= 0 \\
(\x^TA\u-\u^TA\u)^2-(\x^TA\x+\u^TA\u-2\x^TA\u)(\u^TA\u-1) &= 0 \\
(\x^TA\u)^2-2(\x^TA\u)-(\x^TA\x)(\u^TA\u)+\x^TA\x+\u^TA\u &= 0 \\
(\x^TA\u-1)^2-(\x^TA\x-1)(\u^TA\u-1) &= 0.
\end{align*}
In summary, the tangent cone from $\u$ to the surface $\x^TA\x=1$ has the equation
\begin{equation*}
\boxed{(\x^TA\u-1)^2=(\x^TA\x-1)(\u^TA\u-1).}
\end{equation*}
However, since this is a quadric cone with apex at $\u$, we prefer to express it in the form $(\x-\u)^TC(\x-\u)=0$ for some real symmetric matrix $C^T=C$. To find the matrix $C$ we expand the equations in terms of $\x$ to get
\begin{align*}
(\x-\u)^TC(\x-\u) &=0\\
\x^T[C]\x+[-2\u^TC]\x+[\u^TC\u]&=0
\end{align*}
and
\begin{align*}
(\x^TA\u-1)^2-(\x^TA\x-1)(\u^TA\u-1)&=0\\
(\x^TA\u)^2-2(\x^TA\u)-(\x^TA\x)(\u^TA\u)+\x^TA\x+\u^TA\u &= 0 \\
(\x^TA\u)(\x^TA\u)^T-2(\x^TA\u)+(\x^TA\x)(1-\u^TA\u)+\u^TA\u &= 0 \\
\x^T(A\u\u^TA)\x+\x^T\left[(1-\u^TA\u)A\right]\x-2(\x^TA\u)+\u^TA\u &= 0 \\
\x^T\left[A\u\u^TA+(1-\u^TA\u)A\right]\x+[-2\u^TA]\x+[\u^TA\u] &= 0.
\end{align*}
By comparing the leading quadratic forms we see that the cone matrix $C$ satisfies
\begin{equation*}
C=A\u\u^TA+(1-\u^TA\u)A,
\end{equation*}
at least up to a scalar multiple. For this choice of $C$ we easily verify that $-2\u^TC=-2\u^TA$ and $\u^TC\u=\u^TA\u$, so the linear and constant terms also agree.

In summary, we find that the tangent cone from the point $\u$ to the quadric surface $\x^TA\x=1$ has the equation
\begin{equation*}
\boxed{(\x-\u)^T \left[ A\u\u^TA+(1-\u^TA\u)A \right] (\x-\u)=0.}
\end{equation*}

\section{Confocal Quadric Surfaces}
\label{sec:confocal}

The specific surfaces that we need are called {\em confocal quadric surfaces}. This is the key idea that I learned from Hilbert and Cohn-Vossen.

\bigskip

{\bf Definition.} Fix distinct real numbers $a,b,c\in\RR$. Then for any real parameter $k\in\RR$ not in the set $\{a,b,c\}$ we consider the quadric surface
\begin{equation*}
\x^TA_k\x=\frac{x^2}{a-k}+\frac{y^2}{b-k}+\frac{z^2}{c-k}=1,
\end{equation*}
where the matrix $A_k$ is defined by
\begin{equation*}
A_k:=\begin{pmatrix} 1/(a-k) & & \\ & 1/(b-k) & \\ && 1/(c-k) \end{pmatrix}.
\end{equation*}
The surfaces $\x^TA_k\x=1$ for various $k$ describe a {\em confocal family of quadrics}. \hfill ///

\bigskip

For example, let us assume that $c<b<a$. Then the surface $\x^TA_k\x=1$ is
\begin{itemize}
\item an ellipsoid when $k<c$,
\item a hyperboloid of one sheet when $c<k<b$,
\item a hyperboloid of two sheets when $b<k<a$,
\item imaginary when $a<k$.
\end{itemize}
As $k$ approaches one of the critical values $\{a,b,c\}$, the surface degenerates to a space curve in one of the three principal planes of the family. For example, as $k\to c$ we must have $z\to 0$ and hence the surface $\x^TA_k\x=1$ degenerates to the space curve defined by
\begin{equation*}
\frac{x^2}{a-c}+\frac{y^2}{b-c}=1 \quad\text{and}\quad z=0.
\end{equation*}
Note that this is an ellipse in the $x,y$-plane, called the {\em focal ellipse} of the system. For values of $k$ near $c$ we obtain either a very thin ellipsoid on the inside of the ellipse ($k<c$) or a very thin hyperboloid of one sheet on the outside of the ellipse ($c<k$). Similarly, as $k\to b$ or $k\to a$ the surface $\x^TA_k\x=1$ degenerates to the space curve
\begin{equation*}
\frac{x^2}{a-b}+\frac{z^2}{c-b}=1 \quad\text{and}\quad y=0
\end{equation*}
or
\begin{equation*}
\frac{y^2}{b-a}+\frac{z^2}{c-a}=1 \quad\text{and}\quad x=0,
\end{equation*}
respectively. The first of these is a hyperbola in the $x,z$-plane, called the {\em focal hyperbola} of the system. For $k$ on either side of $b$ we obtain very thin hyperboloid of one sheet ($k<b$) or a very thin hyperboloid of two sheets ($b<k$). As $k\to a$ from the left, the two sheets of the hyperboloid converge to the $y,z$-plane from opposite sides and then disappear. The focal curve in the $y,z$-plane is imaginary.

\bigskip

[Remark: The focal {\bf curves} of a confocal family of quadric {\bf surfaces} generalize the focal {\bf points} of a confocal family of quadric {\bf curves}. To see this, one should fix two real numbers $a\neq b$ and consider the family of quadric curves
\begin{equation*}
\frac{x^2}{a-k}+\frac{y^2}{b-k}=1.
\end{equation*}
For more on this topic see Chapter 1 of Hilbert and Cohn-Vossen.]

\bigskip

Here is the ``fundamental theorem" of the subject. We will prove the theorem for surfaces but it should be clear how to generalize the theorem to confocal quadrics in any number of dimensions.

\bigskip

{\bf Fundamental Theorem of Confocal Quadrics.} Fix real numbers $a,b,c\in\RR$ satisfying $c<b<a$ and consider the confocal family $\x^TA_k\x=1$ of quadric surfaces. For each point $\u=(u,v,w)\in\RR^3$ that is not on a principal plane (i.e., with $u,v,w\neq 0$) there exist exactly three surfaces of the family passing through $\u$. These surfaces correspond to parameters $k_1,k_2,k_3\in\RR$ satisfying
\begin{equation*}
k_1<c<k_2<b<k_3<a,
\end{equation*}
hence there is one surface of each topological type. The tangent planes to the three surfaces at $\u$ are mutually perpendicular, and the {\em Cartesian coordinates} $(u,v,w)$ are related to the {\em confocal coordinates} $(k_1,k_2,k_3)$ as follows:
\begin{align*}
u^2 &= (a-k_1)(a-k_2)(a-k_3) / (b-a)(c-a) \\
v^2 &= (b-k_1)(b-k_2)(b-k_3) / (a-b)(c-b) \\
w^2 &= (c-k_1)(c-k_2)(c-k_3) / (a-c)(b-c).
\end{align*}
\hfill ///

\bigskip

Before proving this we must derive the equation of the tangent plane at a given point on a quadric surface. So consider a general quadric surface $\x^TA\x=1$ and let $\u$ be any point on the surface, so that $\u^TA\u=1$. Observe that
\begin{equation*}
\boxed{\x^TA\u=1}
\end{equation*}
is the equation of some plane passing through $\u$. I claim that this is the tangent plane.

To see why, we will show that any line $\u+t\v$ contained in the plane $\x^TA\u=1$ has at least double contact with the surface at $\u$. Indeed, since the line is contained in the plane we have for all $t$ that
\begin{align*}
(\u+t\v)^TA\u &= 1 \\
\u^TA\u+t\,\v^TA\u &= 1\\
1+t\,\v^TA\u &= 1\\
t\,\v^TA\u &=0,
\end{align*}
and it follows that $\v^TA\u=0$. Then the intersection of the line with the surface is determined by the following equation in $t$:
\begin{align*}
(\u+t\v)^TA(\u+t\v) &= 1 \\
\u^TA\u+2t\,\v^TA\u+t^2\,\v^TA\v &= 1\\
1+0+t^2\,\v^TA\v &= 1 \\
t^2\,\v^TA\v &= 0.
\end{align*}
If $\v^TA\v=0$ then we see that the line is completely contained in the surface, and if $\v^TA\v\neq 0$ then we see that $t=0$ is a double root as desired. In summary, we find that the plane $\x^TA\u=1$ is tangent to the surface $\x^TA\x=1$ at $\x=\u$.

\bigskip

{\bf Proof of the Fundamental Theorem.} We have assumed that $c<b<a$. Let $\u=(u,v,w)$ be any point that is not on a principal plane (i.e., such that $u,v,w\neq 0$). We are looking for values of $k$ such that $\u^TA_k\u=1$. In other words, we want
\begin{align*}
\begin{pmatrix} u & v & w\end{pmatrix} \begin{pmatrix} 1/(a-k) && \\ &1/(b-k)& \\ && 1/(c-k) \end{pmatrix} \begin{pmatrix} u\\ v\\ w\end{pmatrix} &= 1\\
\frac{u^2}{a-k}+\frac{v^2}{b-k}+\frac{w^2}{c-k} &=1 \\
(b-k)(c-k)u^2+(a-k)(c-k)v^2+(a-k)(b-k)w^2 &= (a-k)(b-k)(c-k).
\end{align*}
Therefore we will define the polynomial
\begin{equation*}
\boxed{\varphi_\u(k):=(b-k)(c-k)u^2+(a-k)(c-k)v^2+(a-k)(b-k)w^2- (a-k)(b-k)(c-k).}
\end{equation*}
We observe that $\varphi_u(k)$ is a cubic polynomial in $k$ with leading coefficient $1$ and satisfying
\begin{align*}
\varphi_\u(c)&=(a-c)(b-c)w^2>0,\\
\varphi_\u(b)&=(a-b)(c-b)v^2<0,\\
\varphi_\u(a)&=(b-a)(c-a)u^2>0.
\end{align*}
It follows that $\varphi_u(k)$ has three distinct real roots $k_1,k_2,k_3$ satisfying
\begin{equation*}
k_1<c<k_2<b<k_3<a.
\end{equation*}
By previous remarks we know that the surfaces corresponding to $k_1,k_2,k_3$ are an ellipsoid, a hyperboloid of one sheet and a hyperboloid of two sheets, respectively. Thus we have found one surface of each topological type passing through $\u$.

By the remarks after the statement of the theorem, the tangent planes to these three surfaces at $\u$ have the equations
\begin{equation*}
\x^TA_{k_1}\u=1, \quad \x^TA_{k_2}\u=1 \quad\text{and}\quad \x^TA_{k_3}\u=1,
\end{equation*}
respectively. To show that these planes are mutually perpendicular it is enough to show that the normal vectors $A_{k_1}\u$, $A_{k_2}\u$ and $A_{k_3}\u$ are mutually perpendicular, and for this we will use a clever trick. Observe that for any real numbers $k\neq\ell$ with $k,\ell\not\in\{a,b,c\}$ we have the following ``partial fractions" identity:
\begin{equation*}
\boxed{A_kA_\ell=\frac{A_k-A_\ell}{k-\ell}.}
\end{equation*}
Now if $k\neq \ell$ are two elements of the set $\{k_1,k_2,k_3\}$ then by definition we have $\u^TA_k\u=1$ and $\u^TA_\ell\u=1$, and it follows that the vectors $A_k\u$ and $A_\ell\u$ are perpendicular:
\begin{equation*}
(A_k\u)^T(A_\ell\u)=\u^T(A_kA_\ell)\u=\frac{\u^T(A_k-A_\ell)\u}{k-\ell}=\frac{\u^TA_k\u-\u^TA_\ell\u}{k-\ell}=\frac{1-1}{k-\ell}=0.
\end{equation*}
It only remains to solve for the Cartesian coordinates $u,v,w$ in terms of the confocal coordinates $k_1,k_2,k_3$. Since the cubic polynomial $\varphi_\u(k)$ has leading coefficient $1$ and distinct roots $k_1,k_2,k_3$ we must have
\begin{equation*}
\varphi_\u(k)=(k-k_1)(k-k_2)(k-k_3).
\end{equation*}
Then by substituting $k=a,b,c$ into $\varphi_\u(k)$ we obtain
\begin{equation*}
\begin{array}{rcccl}
(b-a)(c-a)u^2 &=& \varphi_\u(a) &=& (a-k_1)(a-k_2)(a-k_3),\\
(a-b)(c-b)v^2 &=& \varphi_\u(b) &=& (b-k_1)(b-k_2)(b-k_3),\\
(a-c)(b-c)w^2 &=& \varphi_\u(c) &=& (c-k_1)(c-k_2)(c-k_3),
\end{array}
\end{equation*}
as desired. \hfill $\qed$

\section{Symmetries of a Tangent Cone}
\label{sec:symmetries}

The final ingredient that we need is the result stated in the footnote of Hilbert and Cohn-Vossen that I alluded to in the Introduction. This fact was apparently well known in the early twentieth century\footnote{Hilbert and Cohn-Vossen refer to the work of Otto Staude but they don't provide a reference. It seems that Staude's work is summarized in his textbook (1896).} but I have not been able to a find a proof in the literature. The selection of topics in this note was chosen to make the proof as slick as possible. Here is the footnote.

\bigskip

{\bf Footnote 4 from Hilbert and Cohn-Vossen (pg. 24).} {\em The following is another property of the confocal system, which, incidentally, includes the property just mentioned\,\footnote{Their ``property just mentioned" was a geometric description of our Main Theorem.} as a limiting case: The planes of symmetry of the tangent cone from any point $\u$ in space to any surface of the system which does not enclose $\u$ are the tangent planes at $\u$ to the three surfaces of the system that pass through $\u$.}

\bigskip

To express this in our language we fix distinct real numbers $c<b<a$ and consider the confocal system of quadric surfaces $\x^TA_k\x=1$ as in Section \ref{sec:confocal}. For a generic point in space $\u\in\RR^3$ (i.e., not on the principal planes) recall from the Fundamental Theorem that there exist three confocal surfaces through $\u$ corresponding to some ``confocal parameters" $k_1,k_2,k_3\in\RR$ satisfying
\begin{equation*}
k_1<c<k_2<b<k_3<c.
\end{equation*}
Furthermore, recall that tangent planes to the three confocal surfaces at $\u$ are mutually perpendicular and are given by the equations
\begin{equation*}
\x^TA_{k_1}\u=1,\quad \x^TA_{k_2}\u=1 \quad\text{and}\quad \x^TA_{k_3}\u=1.
\end{equation*}
Here is the generic case of Hilbert and Cohn-Vossen's statement.

\bigskip

{\bf Theorem.} Let $\u\in\RR^3$ be a generic point (i.e., with nonzero coordinates). Fix any generic confocal parameter $\ell\in\RR\setminus\{a,b,c,k_1,k_2,k_3\}$ and consider the tangent cone from the point $\u$ to the confocal surface $\x^TA_\ell\x=1$. We know from Section \ref{sec:tangentcone} that this cone has equation $(\x-\u)^TK_{\u,\ell}(\x-\u)=0$ where the cone matrix is defined by
\begin{equation*}
\boxed{K_{\u,\ell}:=A_\ell\u\u^T A_\ell+(1-\u A_\ell\u^T)A_\ell.}
\end{equation*}
I claim that this matrix has eigenvectors $A_{k_1}\u$, $A_{k_2}\u$ and $A_{k_3}\u$, with corresponding eigenvalues
\begin{equation*}
\begin{array}{rcl}
\lambda_1 &=& (\ell-k_2)(\ell-k_3)\,/\, (a-\ell)(b-\ell)(c-\ell), \\
\lambda_2 &=& (\ell-k_1)(\ell-k_3)\,/\, (a-\ell)(b-\ell)(c-\ell), \\
\lambda_3 &=& (\ell-k_1)(\ell-k_2)\,/\, (a-\ell)(b-\ell)(c-\ell). \\
\end{array}
\end{equation*}
Since the parameters $k_1,k_2,k_3$ are distinct we observe that the (nonzero) eigenvalues $\lambda_1,\lambda_2,\lambda_3$ are also distinct. If the eigenvalues are not all of the same sign (i.e., if $k_1<\ell<k_3$) then we conclude that $(\x-\u)^TK_{\u,\ell}(\x-\u)=0$ is a real non-circular cone with planes of symmetry equal to the tangent planes $\x^TA_{k_1}\u=1$, $\x^TA_{k_2}\u=1$ and $\x^TA_{k_3}\u=1$ through $\u$. [Remark: It is interesting that that the planes of symmetry depend only on the point $\u$ and not on the parameter $\ell$.] \hfill ///

\bigskip

{\bf Proof.} Consider any $i\in\{1,2,3\}$. To prove that $A_{k_i}\u$ is an eigenvector of $K_{\u,\ell}$ we will use the partial fractions identity
\begin{equation*}
A_\ell A_{k_i} = \frac{A_\ell - A_{k_i}}{\ell-k_i},
\end{equation*}
which holds because $\ell\neq k_i$ and $\ell,k_i\not\in\{a,b,c\}$. Since $\u$ is on the surface $\x^T A_{k_i}\x=1$ (i.e., $\u^TA_{k_i}\u=1$) we have
\begin{align*}
K_{\u,\ell} A_{k_i}\u &= \left[ A_\ell\u\u^T A_\ell+(1-\u A_\ell\u^T)A_\ell\right] A_{k_i}\u \\
&= A_\ell\u\u^TA_\ell A_{k_i}\u +(1-\u A_\ell \u^T) A_\ell A_{k_i}\u \\
&= A_\ell\u\u^T\left( \frac{A_\ell-A_{k_i}}{\ell-k_i}\right)\u+(1-\u A_\ell \u^T)\left( \frac{A_\ell-A_{k_i}}{\ell-k_i}\right)\u \\
&= A_\ell\u\left( \frac{\u^TA_\ell\u-\u^TA_{k_i}\u}{\ell-k_i}\right)+(1-\u A_\ell \u^T)\left( \frac{A_\ell\u-A_{k_i}\u}{\ell-k_i}\right) \\
&= A_\ell\u\left( \frac{\u^TA_\ell\u-1}{\ell-k_i}\right)+(1-\u A_\ell \u^T)\left( \frac{A_\ell\u-A_{k_i}\u}{\ell-k_i}\right) \\
&= \cancel{A_\ell\u\left( \frac{\u^TA_\ell\u-1}{\ell-k_i}\right)}+\cancel{\left( \frac{1-\u^TA_\ell\u}{\ell-k_i}\right)A_\ell\u}-\left(\frac{1-\u^TA_\ell\u}{\ell-k_i}\right) A_{k_i}\u \\
&=\left(\frac{\u^TA_\ell\u-1}{\ell-k_i}\right) A_{k_i}\u.
\end{align*}
It follows that $A_{k_i}\u$ is an eigenvector of $K_{\u,\ell}$ with eigenvalue $\lambda_i=\left(\u^TA_\ell\u-1\right)/(\ell-k_i)$. To compute the eigenvalue explicitly we recall from from the proof of the Fundamental Theorem that the cubic polynomial $\varphi_\u(k):=(a-k)(b-k)(c-k)\left(\u^TA_k\u-1\right)$ has distinct roots $k_1,k_2,k_3$ and hence
\begin{align*}
(a-k)(b-k)(c-k)\left(\u^T A_k\u-1\right) &= (k-k_1)(k-k_2)(k-k_3)\\
\u^TA_k\u-1 &= (k-k_1)(k-k_2)(k-k_3) / (a-k)(b-k)(c-k).
\end{align*}
By substituting $k=\ell$ we obtain the eigenvalue
\begin{equation*}
\lambda_i = \frac{\u^TA_\ell\u-1}{\ell-k_i}=\frac{(\ell-k_1)(\ell-k_2)(\ell-k_3)}{(a-\ell)(b-\ell)(c-\ell)(\ell-k_i)},\end{equation*}
which agrees with the claimed formulas for $\lambda_1,\lambda_2,\lambda_3$ when $i=1,2,3$. \hfill $\qed$

\bigskip

To complete the proof of the Main Theorem it only remains to examine what happens when the point $\u\in\RR^3$ approaches a principal plane and when the confocal parameter $\ell\in\RR$ approaches one of the critical values $\{a,b,c\}$.

\section{Proof of the Main Theorem}
\label{sec:mainthm}

In the previous section we computed the symmetries of the tangent cone from a generic point $\u\in\RR^3$ to a generic surface $\x^TA_\ell\x=1$ in the confocal family corresponding to the fixed parameters $c<b<a$. If $k_1<c<k_2<b<k_3<a$ are the confocal coordinates of the point $\u$, we found that the tangent cone is real when $k_1<\ell<k_3$ and that it is never circular.

In this section we will complete the proof of the Main Theorem by allowing the point $\u$ to approach the principal planes and by allowing the parameter $\ell$ to approach one of the critical values $\{a,b,c\}$, i.e., by allowing the quadric surface $\x^TA_\ell \x=1$ to degenerate to one of the {\em focal curves} of the confocal system:
\begin{align*}
x^2/(a-c)+y^2/(b-c)=1 \quad\text{and}\quad z=0,\\
x^2/(a-b)+z^2/(c-b)=1 \quad\text{and}\quad y=0,\\
y^2/(b-a)+z^2/(c-a)=1 \quad\text{and}\quad x=0.
\end{align*}
Since $c<b<a$ we see that the first of these curves is an ellipse in the $x,y$-plane, the second is a hyperbola in the $x,z$-plane, and the third is an imaginary curve in the $y,z$-plane.

\bigskip

{\bf Letting the surface $\x^TA_\ell\x=1$ degenerate to a focal curve.}

Let $\u\in\RR^3$ be a fixed point not on a principal plane, with confocal coordinates $k_1<c<k_2<b<k_3<a$. Let $(\x-\u)^T K_{\u,\ell} (\x-\u)=0$ be the tangent cone from the point $\u$
to the surface $\x^TA_\ell\x=1$. From the previous section we know that the matrix $K_{\u,\ell}$ has eigenvectors $A_{k_1}\u$, $A_{k_2}\u$ and $A_{k_3}\u$, with corresponding eigenvalues
\begin{equation*}
\begin{array}{rcl}
\lambda_1 &=& (\ell-k_2)(\ell-k_3)\,/\, (a-\ell)(b-\ell)(c-\ell), \\
\lambda_2 &=& (\ell-k_1)(\ell-k_3)\,/\, (a-\ell)(b-\ell)(c-\ell), \\
\lambda_3 &=& (\ell-k_1)(\ell-k_2)\,/\, (a-\ell)(b-\ell)(c-\ell). \\
\end{array}
\end{equation*}
Note that the eigen{\bf vectors} are independent of the parameter $\ell$. The eigen{\bf values} become undefined as $\ell$ approaches one of the critical values $\{a,b,c\}$, however this is easy to fix.

To see what happens as $\ell\to c$ (i.e., as the surface $\x^T A_\ell \x=1$ degenerates to the focal ellipse) we observe that the matrix $(c-\ell)K_{\u,\ell}$ has the same eigenvectors as $K_{\u,\ell}$ but with eigenvalues
\begin{equation*}
\begin{array}{rcl}
\lambda_1 &=& (\ell-k_2)(\ell-k_3)\,/\, (a-\ell)(b-\ell), \\
\lambda_2 &=& (\ell-k_1)(\ell-k_3)\,/\, (a-\ell)(b-\ell), \\
\lambda_3 &=& (\ell-k_1)(\ell-k_2)\,/\, (a-\ell)(b-\ell). \\
\end{array}
\end{equation*}
Since these eigenvalues are well-defined when $\ell\to c$, we conclude that the matrix $K_{\u,c}:=\lim_{\ell\to c} (c-\ell)K_{\u,\ell}$ exists\footnote{Unfortunately it seems that this matrix does not have a nice closed formula.} and is uniquely determined by having eigenvectors $A_{k_1}\u, A_{k_2}\u, A_{k_3}\u$ with corresponding eigenvalues
\begin{equation*}
\begin{array}{rcl}
\lambda_1 &=& (c-k_2)(c-k_3)\,/\, (a-c)(b-c), \\
\lambda_2 &=& (c-k_1)(c-k_3)\,/\, (a-c)(b-c), \\
\lambda_3 &=& (c-k_1)(c-k_2)\,/\, (a-c)(b-c). \\
\end{array}
\end{equation*}
Now the equation $(\x-\u)^T K_{\u,c}(\x-\u)=0$ defines the cone from the point $\u$ to the focal ellipse. Since $k_1<c<k_2<k_3$ we find that the eigenvalues of $K_{\u,c}$ satisfy
\begin{equation*}
\lambda_1>0>\lambda_2>\lambda_3.
\end{equation*}
It follows that the cone from a generic point $\u$ to the focal ellipse is real and non-circular. Below we will see what happens when $\u$ is non-generic.

Similarly we can define the matrices $K_{\u,b}:=\lim_{\ell\to b}(b-\ell)K_{\u,\ell}$ and $K_{\u,a}:=\lim_{\ell\to a}(a-\ell)K_{\u,\ell}$, which both have the same eigenvectors $A_{k_1}\u, A_{k_2}\u$ and $A_{k_3}\u$. The eigenvalues of $K_{\u,b}$ are
\begin{equation*}
\begin{array}{rcl}
\lambda_1 &=& (b-k_2)(b-k_3)\,/\, (a-b)(c-b), \\
\lambda_2 &=& (b-k_1)(b-k_3)\,/\, (a-b)(c-b), \\
\lambda_3 &=& (b-k_1)(b-k_2)\,/\, (a-b)(c-b), \\
\end{array}
\end{equation*}
which satisfy
\begin{equation*}
\lambda_3<0<\lambda_1<\lambda_2.
\end{equation*}
It follows that the cone $(\x-\u)^T K_{\u,b}(\x-\u)=0$ from a generic point $\u$ to the focal hyperbola is real and non-circular.

Finally we observe that the eigenvalues of $K_{\u,a}$ are
\begin{equation*}
\begin{array}{rcl}
\lambda_1 &=& (a-k_2)(a-k_3)\,/\, (b-a)(c-a), \\
\lambda_2 &=& (a-k_1)(a-k_3)\,/\, (b-a)(c-a), \\
\lambda_3 &=& (a-k_1)(a-k_2)\,/\, (b-a)(c-a), \\
\end{array}
\end{equation*}
which satisfy
\begin{equation*}
0<\lambda_1<\lambda_2<\lambda_3.
\end{equation*}
Thus the cone $(\x-\u)^T K_{\u,a} (\x-\u)=0$ from a generic point $\u$ to the imaginary focal curve is imaginary, as expected.

\bigskip

{\bf Letting the point $\u$ approach a principal plane.}

Fix real numbers $c<b<a$ as before and recall from the Fundamental Theorem that each generic point $\u=(u,v,w)\in\RR^3$ is contained in three (mutually perpendicular) confocal quadric surfaces corresponding to some parameters
\begin{equation*}
k_1<c<k_2<b<k_3<a.
\end{equation*}
Conversely, any real numbers $k_1,k_2,k_3$ satisfying these inequalities correspond to three confocal surfaces that intersect (perpendicularly) at the eight points $\u=(u,v,w)$ defined by
\begin{align*}
u^2 &= (a-k_1)(a-k_2)(a-k_3)\,/\, (b-a)(c-a) \\
v^2 &= (b-k_1)(b-k_2)(b-k_3)\,/\, (a-b)(c-b) \\
w^2 &= (c-k_1)(c-k_2)(c-k_3)\,/\, (a-c)(b-c).
\end{align*}
Thus we observe that the point $\u$ approaches the three principal planes precisely when the parameters $k_1,k_2,k_3$ approach the critical values $c,b,a$ from the left:
\begin{equation*}
\begin{array}{ccc}
u\to 0 &\Leftrightarrow& k_1\to c \text{ from the left}, \\
v\to 0 &\Leftrightarrow& k_2\to b \text{ from the left}, \\
w\to 0 &\Leftrightarrow& k_3\to a \text{ from the left}. \\
\end{array}
\end{equation*}
As long as the values $k_1,k_2,k_3$ remain distinct we find that the matrix $K_{\u,\ell}$ of the cone from $\u$ to any surface $\x^T A_\ell \x=1$ (including the degenerate cases when $\ell\in\{a,b,c\}$) has distinct eigenvalues $\lambda_1,\lambda_2,\lambda_3$, and thus it remains non-circular.

Under what conditions do we get a circular cone? In other words, under what conditions does the cone matrix $K_{\u,\ell}$ have a repeated eigenvalue? As the confocal parameters $k_1,k_2,k_3$ move around, we observe from the explicit formulas for the eigenvalues $\lambda_1,\lambda_2,\lambda_3$ that
\begin{equation*}
\lim (\lambda_i-\lambda_j)=0 \quad\Leftrightarrow\quad \lim (k_i-k_j)=0.
\end{equation*}
In other words, the cone from $\u$ to a confocal surface or focal curve becomes circular precisely when two of the confocal parameters $k_1,k_2,k_3$ approach each other. Since the parameters of a generic point satisfy
\begin{equation*}
k_1<c<k_2<b<k_3<a,
\end{equation*}
we see that it is {\bf impossible} for $k_1$ and $k_3$ to approach each other. Thus we have two cases:

\bigskip

{\bf Case 1: The point $\u$ approaches the focal hyperbola.} As $k_2\to b\leftarrow k_3$ we find in the limit that the coordinates of the point $\u=(u,v,w)$ satisfy
\begin{equation*}
\begin{array}{rclcl}
u^2 &=& (a-k_1)(a-b)(a-b)\,/\, (b-a)(c-a) &=& (a-b)(a-k_1)\,/\, (a-c), \\
v^2 &=& (b-k_1)(b-b)(b-b)\,/\, (a-b)(c-b) &=& 0, \\
w^2 &=& (c-k_1)(c-b)(c-b)\,/\, (a-c)(b-c) &=& -(c-b)(c-k_1)\,/\, (a-c). \\
\end{array}
\end{equation*}
This implies that
\begin{equation*}
u^2/(a-b)+w^2/(c-b)=1 \quad\text{and}\quad v=0,
\end{equation*}
which tells us that $\u$ is on the focal hyperbola of the system. At the same time, the eigenvalues of the cone from $\u$ to the surface $\x^T A_\ell \x=1$ approach the values
\begin{equation*}
\begin{array}{rclcl}
\lambda_1 &=& (\ell-b)(\ell-b)\,/\, (a-\ell)(b-\ell)(c-\ell) &=& (b-\ell)\,/\,(a-\ell)(c-\ell), \\
\lambda_2 &=& (\ell-k_1)(\ell-b)\,/\, (a-\ell)(b-\ell)(c-\ell)&=& (k_1-\ell)\,/\,(a-\ell)(c-\ell), \\
\lambda_3 &=& (\ell-k_1)(\ell-b)\,/\, (a-\ell)(b-\ell)(c-\ell)&=& (k_1-\ell)\,/\,(a-\ell)(c-\ell). \\
\end{array}
\end{equation*}
Since $k_1<c<b$ implies $k_1\neq b$ we observe that $\lambda_1\neq\lambda_2=\lambda_3$. The cone is real precisely when $k_1<\ell<b$ and in this case the Theorem at the end of Section \ref{sec:circcone} says that the cone is {\bf circular} with angle of aperture $\theta$ satisfying\begin{equation*}
\cos^2\theta=\frac{\lambda_3}{\lambda_3-\lambda_1}=\frac{k_1-\ell}{k_1-b}.
\end{equation*}
We can view the point $\u$ locally as a function of the parameter $k_1$, which satisfies $k_1<\min\{c,\ell\}$. As $-\infty\leftarrow k_1$ the point $\u$ on the focal hyperbola goes to infinity and we have $\cos^2\theta\to 1$, or $\theta\to 0$. That is, from infinitely far away the surface $\x^T A_\ell \x=1$ looks like a point. If $\ell<c$ (i.e., if the surface $\x^T A_\ell\x=1$ is an ellipsoid) then as $k_1\to \ell$ the point $\u$ approaches an ``umbilic point" on the surface of the ellipsoid and the (circular) tangent cone flattens out into the tangent plane at the umbilic point. In the limiting case $\ell\to c$, the point $\u$ approaches one of the foci $(x,y,z)=(\pm\sqrt{a-b},0,0)$ of the focal ellipse. If $c<\ell$ (i.e., if the surface $\x^T A_\ell \x=1$ is a hyperboloid of one sheet) then as $k_1\to c$ we have $\cos^2\theta\to (c-\ell)/(c-b)$ and the angle of aperture reaches the {\bf maximum} value $\theta=\arccos\sqrt{(c-\ell)/(c-b)}$.

In all cases, the axis of symmetry of the cone is given by the eigenvector
\begin{equation*}
A_{k_1}\u=\left(\frac{u}{a-k_1},\frac{v}{b-k_1},\frac{w}{c-k_1}\right)=\left(\frac{u}{a-k_1},0,\frac{w}{c-k_1}\right),
\end{equation*}
which I claim is {\bf tangent} to the focal hyperbola at the point $\u$. Indeed, let us view the point $\u=(u,v,w)$ locally as a function of the parameter $k_1$. By differentiating the formula for $u^2$ with respect to $k_1$ we obtain
\begin{equation*}
2uu' = (u^2)' = \left[ (a-b)(a-k_1)\,/\, (a-c) \right]' = -(a-b)\,/\, (a-c) = -u^2\,/\, (a-k_1),
\end{equation*}
and hence $u'=(-1/2)\cdot u\,/\, (a-k_1)$. Similarly we see that $v'=(-1/2)\cdot v\,/\, (b-k_1)$ and $w'=(-1/2)\cdot w\,/\, (c-k_1)$, thus the tangent vector to the focal hyperbola at $\u$ is given by
\begin{equation*}
(u',v',w')=\left( \frac{-u}{2(a-k_1)},\frac{-v}{2(b-k_1)},\frac{-w}{2(c-k_1)}\right)=-\frac{1}{2} A_{k_1}\u
\end{equation*}
as desired. \hfill ///

\bigskip

{\bf Case 2: The point $\u$ approaches the focal ellipse.} As $k_1\to c\leftarrow k_2$ we find in the limit that the coordinates of the point $\u=(u,v,w)$ satisfy
\begin{equation*}
\begin{array}{rclcl}
u^2 &=& (a-c)(a-c)(a-k_3)\,/\, (b-a)(c-a) &=& (a-c)(a-k_3)\,/\, (a-b), \\
v^2 &=& (b-c)(b-c)(b-k_3)\,/\, (a-b)(c-b) &=& -(b-c)(b-k_3)\,/\, (a-b), \\
w^2 &=& (c-c)(c-c)(c-k_3)\,/\, (a-c)(b-c) &=& 0.
\end{array}
\end{equation*}
This implies that
\begin{equation*}
u^2/(a-c)+v^2/(b-c)=1 \quad\text{and}\quad w=0,
\end{equation*}
which tells us that $\u$ is on the focal ellipse of the system. At the same time, the eigenvalues of the cone from $\u$ to the surface $\x^T A_\ell \x=1$ approach the values
\begin{equation*}
\begin{array}{rclcl}
\lambda_1 &=& (\ell-c)(\ell-k_3)\,/\, (a-\ell)(b-\ell)(c-\ell) &=& (k_3-\ell)\,/\,(a-\ell)(b-\ell), \\
\lambda_2 &=& (\ell-c)(\ell-k_3)\,/\, (a-\ell)(b-\ell)(c-\ell)&=& (k_3-\ell)\,/\,(a-\ell)(b-\ell), \\
\lambda_3 &=& (\ell-c)(\ell-c)\,/\, (a-\ell)(b-\ell)(c-\ell)&=& (c-\ell)\,/\,(a-\ell)(b-\ell). \\
\end{array}
\end{equation*}
Since $c<b<k_3$ implies $c\neq k_3$ we observe that $\lambda_1=\lambda_2\neq \lambda_3$. The cone is real precisely when $c<\ell<k_3$ and in this case the Theorem at the end of Section \ref{sec:circcone} says that the cone is {\bf circular} with angle of aperture $\theta$ satisfying
\begin{equation*}
\cos^2\theta=\frac{\lambda_1}{\lambda_1-\lambda_3}=\frac{k_3-\ell}{k_3-c}.
\end{equation*}
We can view the point $\u$ locally as a function of the parameter $k_3$, which satisfies $\max\{b,\ell\}<k_3<a$. If $b<\ell$ (i.e., if the surface $\x^TA_\ell\x=1$ is a hyperboloid of two sheets) then as $\ell\leftarrow k_3$ the point $\u$ approaches an ``umbilic point" on the surface and the (circular) tangent cone flattens out into the tangent plane at the umbilic point. If $\ell<b$ (i.e., if the surface $\x^T A_\ell\x=1$ is a hyperboloid of one sheet) then as $b\leftarrow k_3$ we have $\cos^2\theta\to (b-\ell)/(b-c)$ and the angle of aperture reaches the {\bf maximum} value $\theta=\arccos\sqrt{(b-\ell)/(b-c)}$. In the limiting case $\ell\to b$ the point $\u$ approaches one of the foci $(x,y,z)=(\pm\sqrt{a-c},0,0)$ of the focal hyperbola. For any value of $\ell$, the angle of aperture reaches the {\bf minimum} value $\theta=\arccos\sqrt{(a-\ell)/(a-c)}$ when $k_3\to a$, i.e., when $\u$ approaches one of the points $(x,y,z)=(0,\pm\sqrt{b-c},0)$. (Since the point $\u$ is trapped on an ellipse it can't get infinitely far away.)

Finally, the axis of symmetry of the cone is given by the eigenvector $A_{k_3}\u$. Using a similar argument to the previous case we see that this axis is tangent to the focal ellipse at $\u$. \hfill ///

\bigskip

These results hold for any value of $\ell$ as long as the corresponding cone is real. The Main Theorem is just a summary of these results for the cases when $\ell\in\{a,b,c\}$.

\section{Conclusion}
\label{sec:conclusion}

To conclude the note I will answer the three questions from the Introduction in plain language. Let us consider a central and non-degenerate quadric curve in the real $x,y$-plane:
\begin{equation*}
\frac{x^2}{\alpha}+\frac{y^2}{\beta}=1.
\end{equation*}
We assume that the parameters $\alpha,\beta\in\RR$ satisfy $\alpha>\beta$ and are not both negative. Thus our curve is a non-circular ellipse (when $\alpha>\beta>0$) or a non-rectangular hyperbola (when $\alpha>0>\beta$). If we define $(a,b,c):=(\alpha,\beta,0)$ then our curve becomes
\begin{equation*}
\frac{x^2}{a-c}+\frac{y^2}{b-c}=1,
\end{equation*}
which we identify as either the focal ellipse or the focal hyperbola of a certain family of confocal quadric surfaces in $x,y,z$-space.

\bigskip

{\bf Question 1:} From which points in space does our curve look like a circle?

\bigskip

{\bf Answer:} It looks like a circle from points on the the other real focal curve defined by:
\begin{equation*}
\frac{x^2}{a-b}+\frac{z^2}{c-b}=\frac{x^2}{\alpha-\beta}+\frac{z^2}{-\beta}=1 \quad\text{and}\quad y=0.
\end{equation*}
If our curve is an ellipse/hyperbola in the $x,y$-plane then the other focal curve is a hyperbola/ellipse in the $x,z$-plane, passing through the foci of the original curve.

\bigskip

{\bf Question 2:} In which direction should we look to see the circle?

\bigskip

{\bf Answer:} If we are sitting on the focal curve in the $x,z$-plane then we should look in the direction of the tangent line. The focal curve in the $x,y$-plane then looks like a circle centered on this line.

\bigskip

{\bf Question 3:} How big is the circle?

\bigskip

{\bf Answer:} The apparent size of the circle depends on the angle of aperture $\theta$ of the corresponding circular cone.

Suppose our curve is an ellipse ($\alpha>\beta>0$) and that $\u$ lies on the focal hyperbola. As $\u$ approaches one of the foci $(x,y,z)=(\pm\sqrt{\alpha-\beta},0,0)$ of the ellipse, the cone becomes flat and the ellipse looks like an infinitely big circle. As $\u$ goes to infinity the ellipse looks like an infinitesimally small circle.

On the other hand, suppose that our curve is a hyperbola ($\alpha>0>\beta$) and that $\u$ lies on the confocal ellipse. As $\u$ approaches one of the foci $(x,y,z)=(\pm\sqrt{\alpha-\beta},0,0)$ of the hyperbola, the cone becomes flat and the hyperbola looks like an infinitely big circle. As $\u$ approaches one of the points $(x,y,z)=(0,0,\pm\sqrt{-\beta})$ the angle of aperture reaches the {\bf minimum} value $\theta=\arccos\sqrt{\alpha/(\alpha-\beta)}$.

\end{document}